\documentclass{amsart}
\usepackage{amsmath,amsthm,amssymb,IMjournal}
\usepackage{}
\usepackage{amssymb}
\usepackage{graphicx}
\usepackage{amsmath}
\usepackage{epstopdf}
\usepackage{color}
\usepackage{multirow}
\numberwithin{equation}{section}

\begin{document}
\newtheorem{The}{Theorem}[section]

\numberwithin{equation}{section}

\title{Lower bounds for eigenvalues of the Steklov eigenvalue problem with variable coefficients}

\author{\|Yu |Zhang|, \|Hai |Bi|,
        \|Yidu |Yang|}

%\rec {January 27, 2006}

\dedicatory{Cordially dedicated to ...}

\abstract
%%%%%%
%%%7 to 12 lines optimum; references in the abstract should be given in full form, e.g. V. Nov\'ak, M. Novotn\'y: Linear extensions of orderings. Czech. Math. J. 50 (2000), 853--864.
In this paper, using new correction to the Crouzeix-Raviart finite element eigenvalue approximations, we obtain lower eigenvalue bounds for the Steklov eigenvalue problem with variable coefficients on $d$-dimensional domains ($d=2, 3$). In addition, we prove that the corrected eigenvalues asymptotically converge to the exact ones from below whether the eigenfunctions are singular or smooth and whether the eigenvalues are large enough or not. Further, we prove that the corrected eigenvalues still maintain the same convergence order as that of uncorrected eigenvalues. Finally, numerical experiments validate our theoretical results.
\endabstract



\thanks
  Project supported by the Young Scientific and Technical Talents Development of Education Department of Guizhou Province (KY [2018]153), the National Natural Science Foundation of China (Grant No. 11561014 and No. 11761022).
\endthanks

\section{Introduction}\label{sec1}
\indent It is an important topic to obtain upper and lower bounds for eigenvalues. As we all know, thanks to the minimum-maximum principle, it is easy to obtain guaranteed upper bounds of eigenvalues by conforming finite element methods (FEMs). Naturally, attentions have been paid to finding lower bounds of eigenvalues by nonconforming finite elements, such as the rotated bilinear ($Q^{rot}_1$) finite element \cite{lin2008,lin2010,hu2013,liq}, the extension of $Q^{rot}_1$ finite element \cite{lin2008,li2008,liq,hu2014}, the enriched Crouzeix-Raviart (ECR) finite element \cite{hu2013,hu2014,liq,lin2013,luo2012,xie2018}, the Wilson finite element \cite{lin2008,zhang2007}, the Morley element \cite{carstensen2014,hu2014,yang2012}, etc. Especially, a lot of work has been done on the lower bounds for eigenvalues based on the Crouzeix-Raviart (CR) finite element approximations (see \cite{armentano,lin2010,yangzhang,hu2014,yang2015,liq,yangl2009,carstensen1,liux,hu2016,xie2018,youc} and therein).\\
\indent In this paper, we will consider lower bounds for eigenvalues of the following Steklov eigenvalue problem with variable coefficients
  \begin{equation}\label{s2.1}
         \left\{
         \begin{aligned}
         &-div(\alpha\nabla u)+\beta u=0,~~~~~~\text{in}~~ \Omega,\\
         &\alpha\frac{\partial u}{\partial \nu}=\lambda u,~~~~~~~~~~~~~~~~~~~~~~\text{on}~~ \partial\Omega,
         %-\triangle u+u&=0~~~~~~~~\text{in}~~ \Omega, \\
         %\frac{\partial u}{\partial \nu}&=\lambda u~~~~~~\text{on}~~ \partial\Omega,
         \end{aligned}
         \right.
 \end{equation}
 where $\Omega\subset \mathbb{R}^d (d=2, 3)$ is a bounded polygonal domain and $\frac{\partial u}{\partial \nu}$ is the outward normal derivative on $\partial \Omega$. Symbols $\nabla$ and $div$ denote the divergence and the gradient operators,
respectively. $\beta=\beta(x)\in L^\infty(\Omega)$ has positive lower bound, $\alpha=\alpha(x)\in W^{1,\infty}(\Omega)$
  and $\alpha_0\leq \alpha(x)$ for a given constant $\alpha_0>0$.\\
 \indent
   Among the above references, \cite{liq,yangl2009,youc,lin2013} discuss lower eigenvalue bounds for the Steklov eigenvalue problem with constant coefficients. \cite{liq} proves that the CR finite element produces asymptotic lower bounds for eigenvalue in the case of singular eigenfunction. And \cite{liq,yangl2009} also prove that the property of lower bounds in the case of nonsingular eigenfunction but under an additional condition that the eigenvalue is large enough. \cite{youc} obtains guaranteed lower bounds for eigenvalues by correcting the CR finite element eigenvalues approximations, but convergence order of the corrected eigenvalues cannot achieve that of the uncorrected eigenvalues. \cite{lin2013} considers the lower bounds for eigenvalues of the Steklov eigenvalue problem by the ECR finite element (see \cite{hu2014,luo2012}). \\
\indent  Based on the above work, we further discuss asymptotic lower bounds of eigenvalues for the Steklov eigenvalue problem with variable coefficients. We introduce a new correction formula (3.5) to the CR finite element eigenvalues approximations $\lambda_h$ and obtain the corrected eigenvalues $\lambda^c_h$.
Our work has the following features:\\
  \noindent (1)  We prove the following conclusion in Theorem 3.1 (when mesh diameter $h$ is sufficiently small)
  \begin{equation*}
   \lambda\geq\lambda^c_h,
  \end{equation*}
  which shows that the corrected eigenvalues are lower bounds of the exact ones whether the eigenfunctions are singular or smooth and whether the eigenvalues are large enough or not (see Section 3 for details). \\
\noindent (2) The result in Theorem 3.2 implies that the corrected eigenvalues converge to the exact ones without the loss of convergence order, i.e.,
convergence order of the corrected eigenvalues is still the same as that of the uncorrected eigenvalues. \\%Our work is important and meaningful.\\
\noindent (3) For $d$-dimensional domains ($d=2, 3$), we implement numerical experiments in Section 4. Numerical results coincide with the theoretical analysis. We are particularly pleased that the correction takes very little time.\\
\indent It should be pointed out that the correction method and theoretical analysis in
this paper are also valid for the ECR finite element (see Remark 3.1 in Section 3).\\
\indent As for the basic theory of finite element and spectral approximation, we refer
to \cite{babuska1991,boffi2010,oden2012,brenner2002}. Throughout this paper, $C$ denotes a generic positive constant independent on mesh size, which may not be the same at each occurrence.
\section{Preliminary}
\indent
Let $H^s(\Omega)$ denote the Sobolev space with real order $s$ on $\Omega$. Let $\|\cdot\|_{s,\Omega}$ and $|\cdot|_{s,\Omega}$ be the norm and seminorm on $H^s(\Omega)$, respectively. $H^0(\Omega)=L^2(\Omega)$. $H^s(\partial\Omega)$ denotes the Sobolev space with real order $s$ on $\partial\Omega$, $\|\cdot\|_{s,\partial\Omega}$ is the norm on $H^s(\partial\Omega)$ and $H^0(\partial\Omega)=L^2(\partial\Omega)$.\\
 \indent The weak form of (\ref{s2.1}) can be written as: find $(\lambda, u)\in \mathbb{R}\times H^1(\Omega)$, $\|u\|_{0, \partial\Omega}=1$ such that
 \begin{equation}\label{s2.3}
   a(u,v)=\lambda b(u,v),~~~~~~\forall v\in H^1(\Omega),
 \end{equation}
 where
 \begin{equation}\label{s2.4}
 a(u,v)=\int\limits_{\Omega}(\alpha\nabla u \cdot\nabla v+\beta uv)dx,\\
 \end{equation}
 \begin{equation}\label{s2.5}
 b(u,v)=\int\limits_{\partial\Omega}uvds.
 \end{equation}
 \indent Let $\pi_h=\{\kappa\}$ be a regular partition of $\Omega$ with the mesh diameter $h=\max \{h_\kappa\}$ where $h_\kappa$ is the diameter of element $\kappa$. $\varepsilon_h$ is the set of $d-1$ dimensional faces of $\pi_h$. We denote by $|\kappa|$ the measure of the element $\kappa$. \\
 \indent We consider the CR finite element space, proposed by Crouzeix and Raviart \cite{crouzeix}, as follows:
 \begin{eqnarray*}
 V_h=\{v\in L^2(\Omega):~ v|_\kappa\in P_1(\kappa),~ v~ \text{is continuous at the barycenters}\\
  \text{of the $d-1$ dimensional faces of}~ \kappa,~ \forall \kappa\in \pi_h\}.
 \end{eqnarray*}
 Define $\|v\|_h=(\sum\limits_{\kappa\in\pi_h}\|v\|^2_{1,\kappa})^\frac{1}{2}$. $\|v\|_h$ is the norm on $V_h$.\\
The CR finite element approximation of (\ref{s2.3}) is to find $(\lambda_h, u_h)\in \mathbb{R}\times V_h$, $\|u_h\|_{0,\partial\Omega}=1$, such that
\begin{equation}\label{s2.6}
  a_h(u_h, v)=\lambda_h b(u_h, v),~~~~~~\forall v\in V_h,
\end{equation}
  where
  \begin{equation}\label{s2.7}
 a_h(u_h,v)=\sum\limits_{\kappa\in\pi_h}\int\limits_{\kappa}(\alpha\nabla u_h \cdot\nabla v+\beta u_hv)dx.
 \end{equation}
From Theorem 4 in \cite{savar1998} and Remark 2.1 in \cite{garau}, we have the following regularity result.\\
\noindent {\bf Regularity}:~~Assume that $\varphi$ is the solution of source problem associated with (\ref{s2.3}). If $f\in L^2(\partial\Omega)$, then $\varphi\in H^{1+r}(\Omega)$ for all $r\in (0, \frac{1}{2})$ and
\begin{equation*}
  \|\varphi\|_{1+r}\leq C\|f\|_{0,\partial\Omega}.
\end{equation*}

\noindent{\bf Lemma 2.1}.~
Let $(\lambda_h,u_h)$ be the $j$th eigenpair of (\ref{s2.6}) and $\lambda$ be the $j$th eigenvalue of (\ref{s2.3}). If $h$ is sufficiently small, there exists $u \in H^{1+r}(\Omega)$ such that
\begin{eqnarray}\label{s2.8}
~~~~~~~~~~~~~~~&~& \|u_h-u\|_h\leq Ch^r,\\ \label{s2.10}
%~~~~~~~~~~~~~~~&~& \|u_h-u\|_{0,\partial\Omega}\leq Ch^{2r},\\
~~~~~~~~~~~~~~~&~& |\lambda_h-\lambda|\leq Ch^{2r},%\\\label{s2.11}
%~~~~~~~~~~~~~~~&~&\|u-u_h\|_{0,\Omega}\leq Ch^r\|u-u_h\|_h,
\\\label{s2.12}
~~~~~~~~~~~~~~~&~&\|u-u_h\|_{0,\partial\Omega}\leq Ch^r\|u-u_h\|_h.
\end{eqnarray}
\noindent{\bf Proof}.~When $\Omega\subset \mathbb{R}^2$, using Theorem 3.1 in \cite{liq}, Theorem 4.6 in \cite{alonso2009} and Theorem 2.2 in \cite{russo2011}, we can deduce (\ref{s2.8}) and (\ref{s2.10}). Referring to Lemma 2.3 in \cite{yang2015} and using Nitsche technique, we can obtain (\ref{s2.12}); when $\Omega\subset \mathbb{R}^3$, using similar arguments to the case of $\Omega\subset \mathbb{R}^2$ (as well as referring to Theorem 4 in \cite{yang}), we can prove that the lemma is valid.~~~$\Box$\\

Define the Crouzeix-Raviart interpolation operator $I_h: H^1(\Omega)\rightarrow V_h$ by
\begin{equation}\label{s2.13}
  \int\limits_eI_hu ds=\int\limits_euds,~~~~\forall e\in \varepsilon_h, u\in H^1(\Omega).
\end{equation}
Note that the interpolation operator $I_h$ has an important orthogonality property (see \emph{equality (2.9)} in \cite{armentano}): for each element $\kappa\in \pi_h$, there holds
\begin{equation}\label{s2.14}
  \int\limits_{\kappa}\nabla (u-I_hu)\cdot \nabla v_hdx=\int\limits_{\partial\kappa}(u-I_hu)\nabla v_h\cdot \nu ds=0,~~~~\forall v_h\in V_h.
\end{equation}

\indent The estimation of constants in Poincar\'{e} and the trace inequalities is a concern of academe (e.g., see \cite{sebestova2014,carsten2012,carstensen2014,chavel1977,li2018,youc,liux} and therein). From Theorem 4.2 in \cite{liux}, we have the following Lemma 2.2.\\
\noindent{\bf Lemma 2.2}.~
For any element $\kappa$, the following conclusion is valid:
\begin{equation}\label{s2.15}
  \|u-I_hu\|_{0,\kappa}\leq C_{h_{\kappa}}|u-I_hu|_{1,\kappa},~~~~\forall u\in H^1(\kappa),
\end{equation}
here
 \begin{itemize}
 \item $C_{h_{\kappa}}=0.1893h_\kappa$ \text{for a triangle element} $\kappa$ \text{in} $\mathbb{R}^2$,
 \item $C_{h_{\kappa}}=0.3804h_\kappa$ \text{for a tetrahedron element} $\kappa$ \text{in} $\mathbb{R}^3$.
 \end{itemize}

\indent Consider any element $\kappa$ with nodes $P_1$, $P_2$, $\cdots$, $P_{d+1}$. The opposite edge/face of point $P_{d+1}$ is denoted by $e$. The measure of $e$ is $|e|$. $H_\kappa$ is the height of element $\kappa$ respect to $e$. It is easy to know $H_\kappa=\frac{d|\kappa|}{|e|}$. Thanks to Lemma 2 of \cite{carstensen2014} and Theorem 3.3 of \cite{youc}, we have the following Lemma 2.3.\\
\noindent{\bf Lemma 2.3}.~
For a given element $\kappa$, there holds:
\begin{equation}\label{s2.16}
  \|u-I_hu\|_{0,e}\leq C_{h_{e}}|u-I_hu|_{1,\kappa},~~~~\forall u\in H^1(\kappa),
\end{equation}
here
 \begin{itemize}
   \item $C_{h_{e}}=0.6711\frac{h_\kappa}{\sqrt{H_\kappa}}$ for a triangle element $\kappa$ in $\mathbb{R}^2$,
    \item $C_{h_{e}}=1.0931\frac{h_\kappa}{\sqrt{H_\kappa}}$ for a tetrahedron element $\kappa$ in $\mathbb{R}^3$.
 \end{itemize}
\noindent{\bf Proof}.~The proof can be found in Theorem 3.3 of \cite{youc}. For convenience of reading, in the case of $d=3$, we write the proof here again.
%\begin{figure}[h!]~~~
%\centering
%\includegraphics[width=12cm]{tetrahedron.eps}
%\caption{\emph{tetrahedron element $\kappa$}}
%\end{figure}
%\noindent {\bf Theorem}. For a given tetrahedron element $\kappa$ as shown in Figure 1, the following error estimate holds for any $u\in H^1(\kappa)$:
%\begin{equation*}
%  \|u-I_hu\|_{0,e}\leq \sqrt{\frac{|e|}{3|\kappa|}(3C^2_{h_\kappa}+2h_\kappa C_{h_\kappa})}|u-I_hu|_{1,\kappa}.
%\end{equation*}
For any $v\in H^1(\kappa)$, from Green formula we have
\begin{equation}\label{s1.1}
 \int\limits_{\kappa}((x_1,x_2,x_3)-P_4)\cdot \nabla(v^2)d\kappa=\int\limits_{\partial\kappa}((x_1,x_2,x_3)-P_4)\cdot \mathbf{n}v^2ds-\int\limits_{\kappa}3v^2d\kappa,
\end{equation}
We deduce
\begin{equation}\label{s1.2}
  ((x_1,x_2,x_3)-P_4)\cdot \mathbf{n}=\left\{
     \begin{array}{ll}
       0, & \text{on faces}~P_1P_2P_4,~P_1P_3P_4,~\text {and}~P_2P_3P_4, \\
       \frac{3|\kappa|}{|e|}, & \text{on face}~P_1P_2P_3.
     \end{array}
   \right.
\end{equation}
Substituting (\ref{s1.2}) into (\ref{s1.1}), we obtain
\begin{eqnarray}\nonumber
  \frac{3|\kappa|}{|e|}\int\limits_{e}v^2ds&=&\int\limits_{\kappa}3v^2d\kappa+\int\limits_{\kappa}((x,y,z)-P_4)\cdot\nabla(v^2)d\kappa\\\nonumber
&\leq&3\int\limits_{\kappa}v^2d\kappa+\int\limits_{\kappa}|(x,y,z)-P_4||\nabla(v^2)|d\kappa\\\nonumber
&\leq&3\int\limits_{\kappa}v^2d\kappa+2h_\kappa\int\limits_{\kappa}|v||\nabla(v)|d\kappa\\\label{s1.3}
&\leq&3\|v\|^2_{0,\kappa}+2h_\kappa\|v\|_{0,\kappa}\|\nabla(v)\|_{0,\kappa},
\end{eqnarray}
Taking $v=u-I_hu$ and applying estimate (\ref{s2.15}) we deduce
\begin{equation*}
  \|u-I_hu\|^2_{0,e}\leq \frac{|e|}{3|\kappa|}(3C^2_{h_\kappa}+2h_\kappa C_{h_\kappa})|u-I_hu|^2_{1,\kappa},
\end{equation*}
which implies that (\ref{s2.16}) is valid when $\Omega\subset \mathbb{R}^3$.~~~$\Box$
\section{The lower bounds property of corrected eigenvalues}
\indent For the problem (\ref{s2.1}), thanks to the minimum-maximum principle, it is easy to obtain guaranteed upper bounds for eigenvalues by conforming finite element methods. From \cite{liq}, we know that CR finite element method gives asymptotic lower
bounds for eigenvalues when the corresponding eigenfunctions are singular or the eigenvalues are large enough. In this section, we introduce a correction for eigenvalues of the problem (\ref{s2.1}) and we will prove that the corrected eigenvalues converge to the exact ones from below whether the corresponding eigenfunctions are smooth or singular and whether the eigenvalues are large enough or not. For preparation, we prove the following inequality (\ref{s3.1}) and Lemma 3.1.\\
 Using (\ref{s2.14}) we have $$\int\limits_{\kappa}\nabla (u-I_hu)\cdot\nabla (u-I_hu)dx=\int\limits_{\kappa}\nabla (u-I_hu)\cdot\nabla (u-u_h)dx\leq|u-I_hu|_{1,\kappa}|u-u_h|_{1,\kappa},$$ then
\begin{equation}\label{s3.1}
  |u-I_hu|_{1,\kappa}\leq |u-u_h|_{1,\kappa}.
\end{equation}
The identity in the following Lemma 3.1 is an equivalent form
of the identity (4.1) in \cite{liq}, which is a generalization of the identities (2.12) in \cite{armentano} and (2.3) in \cite{zhang2007}.\\
\noindent{\bf Lemma 3.1}.~
Let $(\lambda, u)$ and $(\lambda_h, u_h)$  be an eigenpair of (\ref{s2.3}) and (\ref{s2.6}), respectively. Then the following identity is valid:
\begin{eqnarray}\nonumber
  ~~~~~~~~ \lambda-\lambda_h&=&a_h(u-u_h, u-u_h)-\lambda_h b(u-u_h, u-u_h)\\\label{s3.2}
   ~~~~~~~~ &~&+2a_h(u-I_hu, u_h)-2\lambda_hb(u-I_hu,u_h).
\end{eqnarray}
\noindent{\bf Proof}.~
From $\|u\|_{0,\partial\Omega}=1=\|u_{h}\|_{0,\partial\Omega}$, we get
\begin{equation*}
a_{h}(u,u)=\lambda,~~~~a_{h}(u_{h},u_{h})=\lambda_{h}.
\end{equation*}
Therefore
\begin{eqnarray}\label{s3.3}
\lambda-\lambda_{h}&=&a_{h}(u,u)+a_{h}(u_{h},u_{h})-2a_{h}(u_{h},u_{h})\nonumber\\
&=&a_{h}(u,u)+a_{h}(u_{h},u_{h})-2a_{h}(u,u_{h})+2a_{h}(u-u_{h},u_{h})\nonumber\\
&=&a_{h}(u-u_{h},u-u_{h})+2a_{h}(u-u_{h},u_{h}).
\end{eqnarray}
From $b(I_{h}u-u_{h},u_{h})=b(I_{h}u-u,u_{h})+b(u-u_{h},u_{h}-\frac{1}{2}u+\frac{1}{2}u)$, we obtain
\begin{equation*}
\lambda_{h}b(I_{h}u-u_{h},u_{h})
=\lambda_{h}b(I_{h}u-u,u_{h})-\frac{1}{2}\lambda_{h}b(u-u_{h},u-u_{h}),
\end{equation*}
which together with (\ref{s2.6}) yields
\begin{eqnarray}\label{s3.4}
&&a_{h}(u-u_{h},u_{h})=a_{h}(u-I_{h}u,u_{h})+a_{h}(I_{h}u-u_{h},u_{h})\nonumber\\
&~&~~~~~~=a_h(u-I_{h}u,u_{h})+\lambda_{h}b(I_{h}u-u_{h},u_{h})\nonumber\\
&~&~~~~~~=a_h(u-I_{h}u,u_{h})+\lambda_{h}b(I_{h}u-u,u_{h})-\frac{1}{2}\lambda_{h}b(u-u_{h},u-u_{h}).
\end{eqnarray}
Substituting (\ref{s3.4}) into (\ref{s3.3}), we get (\ref{s3.2}).~~~$\Box$\\

\indent Now we give correction formula (\ref{s3.5}). In addition, we will prove that the correction provides asymptotic
lower bounds for eigenvalues of the problem (\ref{s2.3}).\\
\indent Denote by $I_0$ the piecewise constant interpolation operator on $\Omega$.
Let $(\lambda, u)$ be an eigenpair of (\ref{s2.3}) and $(\lambda_h, u_h)$ be the corresponding CR finite element approximations. We introduce the following formula to correct the CR finite element approximations $\lambda_h$:
\begin{equation}\label{s3.5}
\lambda^c_h=\frac{\lambda_h}{1+\frac{\delta}{\lambda_h\alpha_0}\sum\limits_{\kappa\in\pi_h}(\|(\alpha-I_0\alpha)\nabla u_h\|_{0,\kappa}+C_{h_{\kappa}}\|\beta u_h\|_{0,\kappa})^2},
\end{equation}
where $\delta>1$ is any given constant.\\
\indent For the convenience of the next proof, we denote
$$M=\frac{\delta}{\alpha_0}\sum\limits_{\kappa\in\pi_h}(\|(\alpha-I_0\alpha)\nabla u_h\|_{0,\kappa}+C_{h_{\kappa}}\|\beta u_h\|_{0,\kappa})^2,$$
then
\begin{equation}\label{z3.12}
 \lambda^c_h=\frac{\lambda_h}{1+\frac{1}{\lambda_h}M}.
\end{equation}
By the interpolation error estimate, we know
\begin{equation}\label{s3.11}
  \|\alpha-I_0\alpha\|_{0,\infty,\kappa}\leq Ch_\kappa\|\alpha\|_{1,\infty,\kappa}.
\end{equation}
Noting that $C_{h_{\kappa}}=0.1893h_\kappa$, we derive %and $C_{h_{e}}=0.6711\frac{h_\kappa}{\sqrt{H_\kappa}}$, we derive
\begin{equation}\label{s3.12}
 0\leq M\leq Ch^2.
\end{equation}

\indent In practical computation, we can't guarantee that $\lambda_h$ are lower bounds of $\lambda$ if we are not sure the eigenfunctions are singular or the eigenvalues are large enough. Now we will prove the corrected eigenvalues $\lambda^c_h$ are lower bounds of the exact ones whether the eigenfunctions are singular or smooth and whether the eigenvalues are large
enough or not.\\
\noindent{\bf Theorem 3.1}.~
Let $\lambda^c_h$ be a corrected eigenvalue obtained by (\ref{s3.5}). Assuming that the conditions of Lemma 2.1 hold and $\|u-u_h\|_{h}\geq Ch^{1 +\frac{r}{2}}$, then we have the following conclusion:
\begin{equation}\label{s3.6}
\lambda\geq\lambda^c_h.
\end{equation}

\noindent{\bf Proof}.~
We discuss the four terms on the right-hand side of (\ref{s3.2}). Since $\alpha\geq\alpha_0$, we have
\begin{equation}\label{s3.7}
a_h(u-u_h, u-u_h) \geq \sum\limits_{\kappa\in\pi_h}(\alpha_0|\nabla(u-u_h)|^2_{1,\kappa}+\int\limits_{\kappa}\beta(u-u_h)^2dx).
  %a_h(u-u_h, u-u_h) = \sum\limits_{\kappa\in\pi_h}(|u-u_h|^2_{1,\kappa}+\|u-u_h\|^2_{0,\kappa}).
\end{equation}
%\indent (\ref{s2.8}) and (\ref{s2.9}) indicate that, on the right-hand side of (\ref{s3.2}), the second term is \\
From (\ref{s2.14}), we have
\begin{eqnarray*}\nonumber
a_h(I_hu-u, u_h)&=&\sum\limits_{\kappa\in\pi_h}\int\limits_{\kappa}\big((\alpha-I_0\alpha)\nabla(I_hu-u)\cdot\nabla u_h
+I_0\alpha\nabla(I_hu-u)\cdot\nabla u_h\\&~&+\beta(I_hu-u)u_h\big)dx\\
&=&\sum\limits_{\kappa\in\pi_h}\int\limits_{\kappa}((\alpha-I_0\alpha)\nabla(I_hu-u)\cdot\nabla u_h+\beta(I_hu-u)u_h)dx.
\end{eqnarray*}
Applying Cauchy-Schwarz inequality and (\ref{s2.15}) to the above equality, we deduce
\begin{eqnarray*}\nonumber
  a_h(I_hu-u, u_h)%&=&\sum\limits_{\kappa\in\pi_h}\int\limits_{\kappa}(\nabla(I_hu-u)\cdot\nabla u_h+(I_hu-u)u_h)dx\\\nonumber
                  %&=&\sum\limits_{\kappa\in\pi_h}\int\limits_{\kappa}(I_hu-u)u_hdx\\\nonumber
                  &\leq&\sum\limits_{\kappa\in\pi_h}\big(|u-I_hu|_{1,\kappa}\|(\alpha-I_0\alpha)\nabla u_h\|_{0,\kappa}+\|u-I_hu\|_{0,\kappa}\|\beta u_h\|_{0,\kappa}\big)\\
                  %&\leq&\sum\limits_{\kappa\in\pi_h}|u-I_hu|_{1,\kappa}(\|(\alpha-I_0\alpha)\nabla u_h\|_{0,\kappa}+C_{h_{\kappa}}\|\beta u_h\|_{0,\kappa})\\
                  &\leq&\sum\limits_{\kappa\in\pi_h}|u-I_hu|_{1,\kappa}(\|(\alpha-I_0\alpha)\nabla u_h\|_{0,\kappa}+C_{h_{\kappa}}\|\beta u_h\|_{0,\kappa}),
\end{eqnarray*}
which together with Young inequality yields
\begin{equation}\label{s3.8}
 2a_h(I_hu-u, u_h) \leq\frac{\alpha_0}{\delta}\sum\limits_{\kappa\in\pi_h}|u-I_hu|^2_{1,\kappa}+\frac{\delta}{\alpha_0}\sum\limits_{\kappa\in\pi_h}(\|(\alpha-I_0\alpha)\nabla u_h\|_{0,\kappa}+C_{h_{\kappa}}\|\beta u_h\|_{0,\kappa})^2.
\end{equation}
For the later proof, we introduce the piecewise constant interpolation operator $I^b_0$ on $\partial\Omega$.
Using (\ref{s2.13}), Cauchy-Schwarz inequality, (\ref{s2.16}), (\ref{s3.1}), interpolation error estimates and trace inequality, we get
 \begin{eqnarray}\nonumber
 b(u-I_hu,u)&=&\sum\limits_{e\in\varepsilon_h\cap\partial\Omega}\int\limits_{e}((u-I_hu)(u-I^b_0u)+(u-I_hu)I^b_0u)ds\\\nonumber
 &\leq&\sum\limits_{e\in\varepsilon_h\cap\partial\Omega}\|u-I_hu\|_{0,e}\|u-I^b_0u\|_{0,e}\\\nonumber
&\leq& Ch^{\frac{1}{2}+r}\|u\|_{{\frac{1}{2}+r},\partial\Omega}\big(\sum\limits_{\kappa\in\pi_h,e\in\partial\kappa\cap\partial\Omega}C^2_{h_{e}}|u-I_hu|^2_{1,\kappa}\big)^\frac{1}{2}\\\label{zz1}
 &\leq&Ch^{1+r}(\sum\limits_{\kappa\in\pi_h}|u-u_{h}|_{1,\kappa}^2)^{\frac{1}{2}}\|u\|_{1+r}.
\end{eqnarray}
According to $\|u-u_h\|_{h}\geq Ch^{1 +\frac{r}{2}}$ and (\ref{zz1}), we have
 \begin{equation}\label{zz2}
 b(u-I_hu,u)
 \leq Ch^\frac{r}{2}\|u-u_{h}\|_{h}^2.
\end{equation}
From Cauchy-Schwarz inequality, (\ref{s2.16}) and (\ref{s2.12}) we have
\begin{eqnarray}\nonumber
 | b(u-I_hu,u_h-u)|&\leq& \sum\limits_{e\in\varepsilon_h\cap\partial\Omega}\|u-I_hu\|_{0,e}\|u_h-u\|_{0,e}\\\nonumber
 &\leq& Ch^{\frac{1}{2}}(\sum\limits_{\kappa\in\pi_h}|u-I_hu|_{1,\kappa}^{2})^{\frac{1}{2}}h^{r}\|u_h-u\|_{h}\\\label{zz3}
 %&\leq& Ch^{\frac{1+3r}{2}}(\sum\limits_{\kappa\in\pi_h}|u-u_h|_{1,\kappa}^2)^{\frac{1}{2}}\|u\|_{1+r}.
  &\leq& Ch^{\frac{1}{2}+r}\|u_h-u\|^2_{h}.
 \end{eqnarray}
Combining (\ref{zz2}) and (\ref{zz3}), we deduce
\begin{equation}\label{s3.9}
 2\lambda_hb(u-I_hu,u_h)\leq Ch^\frac{r}{2}\|u_h-u\|^2_{h}. %{\color{red}Ch^{\frac{1}{2}+r}(\sum\limits_{\kappa\in\pi_h}|u-u_{h}|_{1,\kappa}^2)^{\frac{1}{2}}\|u\|_{1+r}.}
\end{equation}
From (\ref{s3.2}), (\ref{s3.7}), (\ref{s3.1}), (\ref{s3.8}) and (\ref{s3.9}), we deduce
\begin{eqnarray*}
\lambda-\lambda_h&\geq&(1-\frac{1}{\delta})\alpha_0\sum\limits_{\kappa\in\pi_h}|u-u_h|^2_{1,\kappa}+\sum\limits_{\kappa\in\pi_h}\int\limits_\kappa\beta(u-u_h)^2dx-\lambda_h\|u-u_h\|^2_{0,\partial\Omega}~~~~~~~~~~~\\
&~&-\frac{\delta}{\alpha_0}\sum\limits_{\kappa\in\pi_h}(\|(\alpha-I_0\alpha)\nabla u_h\|_{0,\kappa}+C_{h_{\kappa}}\|\beta u_h\|_{0,\kappa})^2-Ch^\frac{r}{2}\|u_h-u\|^2_{h}.
%-Ch^{\frac{5}{4}}(\sum\limits_{\kappa\in\pi_h}|u-u_{h}|_{1,\kappa}^2)^{\frac{1}{2}}\|u\|_{1+r}.
\end{eqnarray*}
From the definition of $M$, we have
\begin{eqnarray*}
\lambda-\lambda_h&\geq&(1-\frac{1}{\delta})\alpha_0\sum\limits_{\kappa\in\pi_h}|u-u_h|^2_{1,\kappa}+\sum\limits_{\kappa\in\pi_h}\int\limits_\kappa\beta(u-u_h)^2dx-\lambda_h\|u-u_h\|^2_{0,\partial\Omega}\\
&~&-Ch^\frac{r}{2}\|u_h-u\|^2_{h}-\frac{\lambda_h-\lambda}{\lambda_h}M-\frac{\lambda}{\lambda_h}M,
%-Ch^{\frac{5}{4}}(\sum\limits_{\kappa\in\pi_h}|u-u_{h}|_{1,\kappa}^2)^{\frac{1}{2}}\|u\|_{1+r}-\frac{\lambda_h-\lambda}{\lambda_h}M-\frac{\lambda}{\lambda_h}M,
\end{eqnarray*}
which implies that
\begin{eqnarray}\nonumber
(1+\frac{1}{\lambda_h}M)\lambda-\lambda_h&\geq&(1-\frac{1}{\delta})\alpha_0\sum\limits_{\kappa\in\pi_h}|u-u_h|^2_{1,\kappa}+\sum\limits_{\kappa\in\pi_h}\int\limits_\kappa\beta(u-u_h)^2dx\\\label{z3.9}
&~&-\lambda_h\|u-u_h\|^2_{0,\partial\Omega}-Ch^{\frac{1}{2}+r}\|u_h-u\|^2_{h}-\frac{\lambda_h-\lambda}{\lambda_h}M.
%-\lambda_h\|u-u_h\|^2_{0,\partial\Omega}-Ch^{\frac{5}{4}}(\sum\limits_{\kappa\in\pi_h}|u-u_{h}|_{1,\kappa}^2)^{\frac{1}{2}}\|u\|_{1+r}-\frac{\lambda_h-\lambda}{\lambda_h}M.\\
\end{eqnarray}
According to (\ref{s2.12}), it is easy to know that, when $h$ is sufficiently small, the third on the right-hand side of (\ref{z3.9}) are infinitesimals of higher order compared with the sum of the first two terms. From (\ref{s3.12}) and (\ref{s2.10}), we get that the fifth term on the right-hand side of (\ref{z3.9}) is an infinitesimal of higher order compared with the sum of the first two terms. Hence the sign of the right-hand side of (\ref{z3.9}) is determined by summation of the first two terms, i.e., $$(1+\frac{1}{\lambda_h}M)\lambda-\lambda_h\geq 0.$$
From (\ref{z3.12}), we knot that (\ref{s3.6}) is valid. The proof is completed.~~~$\Box$\\

\indent The following theorem shows that $\lambda^c_h$ converge to $\lambda$ and maintain the same convergence order as $\lambda_h$.\\
\noindent{\bf Theorem 3.2}.~
Let $(\lambda, u)$ and $(\lambda_h, u_h)$ be an eigenpair of (\ref{s2.3}) and (\ref{s2.6}), respectively. $\lambda^c_h$ is a corrected eigenvalue obtained by (\ref{s3.5}), then we have
\begin{equation}\label{s3.10}
 \lambda-\lambda^c_h= \lambda-\lambda_h+\frac{\lambda_hM}{\lambda_h+M},
\end{equation}
where $|M|\leq Ch^2$.\\
\noindent{\bf Proof}.~It is easy to deduce the conclusion. From (\ref{z3.12}), we have
\begin{equation*}
  \lambda-\lambda^c_h=\lambda-\lambda_h+\lambda_h-\frac{\lambda_h}{1+\frac{1}{\lambda_h}M}=\lambda-\lambda_h+\frac{\lambda_hM}{\lambda_h+M}.
\end{equation*}
The proof is completed.~~~$\Box$\\

\indent In \cite{liq,lin2013}, it has been obtained that the ECR finite element can produce lower eigenvalue bounds for the Steklov eigenvalue with constant coefficient whether the eigenfunctions are smooth or singular. However, the ECR element cannot produce lower eigenvalue bounds for the Steklov eigenvalue problem with variable coefficients. Therefore we introduce correction to the ECR finite element eigenvalue approximations to obtain lower bounds of eigenvalues.\\
\noindent{\bf Remark 3.1}~(The correction to the ECR finite element eigenvalue approximations).\\
Let $(\lambda_{h},u_{h})$ be approximation eigenpair of (\ref{s2.3}) obtained by ECR element,
\begin{equation*}
\lambda_{h}^c=\frac{\lambda_{h}}{1+\frac{\delta}{\lambda_h\alpha_{0}}\sum\limits_{\kappa\in\pi_{h}}\|(\alpha-I_{o}\alpha)\nabla u_{h}\|_{0,\kappa}^2}.
\end{equation*}
Let
$\beta\in W^{1,\infty}(\Omega)$, then
\begin{eqnarray*}
&& \sum\limits_{\kappa\in\pi_{h}}\int\limits_{\kappa}\beta(u-I_{h}u)u_{h}dx=
  \sum\limits_{\kappa\in\pi_{h}}\int\limits_{\kappa}(u-I_{h}u)(\beta u_{h}-I_{0}(\beta u_{h}))dx\\
&&~~~\leq C  \sum\limits_{\kappa\in\pi_{h}}h_{\kappa}^{2}|u-I_{h}u|_{1,\kappa}\|\beta u_{h}\|_{1,\kappa}.
 \end{eqnarray*}
And using similar argument to Theorems 3.1 and 3.2, we can deduce that when $\|u-u_{h}\|_{h}\geq C h^{1+\frac{r}{2}}$ and $h$ is sufficiently small,
$$\lambda\geq\lambda_{h}^c,$$
and $\lambda_{h}^c$ maintain the same convergence order as $\lambda_{h}$.
\section{Numerical experiments}
\indent In this section, to validate the theoretical results in this paper, we execute correction (\ref{s3.5}) to (\ref{s2.1}) on domain $\Omega$. In computation, we choose $\alpha=\beta=1$. The discrete eigenvalue problems are solved in MATLAB 2018b on an Lenovo ideaPad PC with 1.8GHZ CPU and 8GB RAM. Our program is compiled under the package of iFEM \cite{chen}. The following notations are adopted in tables and figures.
\begin{description}
\item[$h_0$]~~~~The diameter of $\Omega$.
\item[$h$]~~~~~~~~~~~The diameter of meshes.
\item[$\lambda_{j}$]~~~~The $j$th eigenvalue of (\ref{s2.3}).
\item[$\lambda_{j, h}$]~~The $j$th eigenvalue of (\ref{s2.6}) computed by CR finite element.
\item[$\lambda^c_{j,h}$]~~The approximation obtained by correcting $\lambda_{j, h}$.
\item[$t(s)$]~~The CPU time to compute eigenvalues on the finest meshes.
\end{description}
\subsection{Numerical results on $\Omega\subset \mathbb{R}^2$}
\indent When $\Omega\subset \mathbb{R}^2$, we compute on the unit square $(0,1)^2$ ($h_0=\sqrt{2}$), the L-shaped domain $(-1, 1)^2\setminus((0, 1)\times(-1, 0))$ ($h_0=2\sqrt{2}$) and the regular hexagon with side length of 1 ($h_0=2$); for convenience, we simplify the domains as {\bf S}, {\bf L} and {\bf H}, respectively.\\
\indent First, from \cite{liq} we know that CR finite element provides asymptotic lower bounds for eigenvalues of the problem (\ref{s2.1}) when the eigenfunctions are singular or the eigenvalues are large enough. It's worth noting that, on the unit square and the L-shaped domain, $\lambda_{1,h}$ converge to $\lambda_1$ from above in Tables 1 and 5 of \cite{liq}, which imply that the corresponding eigenfunctions are smooth. For the regular hexagon, the same result is obtained. In addition, from Table 2 in \cite{youc}, we know that the guaranteed lower bounds can only achieve convergence order $\mathcal{O}(h)$ even for convex domain. In order to obtain asymptotic lower bounds with the optimal convergence order for the problem (\ref{s2.1}), we use (\ref{s3.5}) to correct $\lambda_{1,h}$. New approximate eigenvalues $\lambda^c_{1,h}$ are listed in Table 1.
  We depict the error curves of $\lambda_{1,h}$, $\lambda^c_{1,h}$ and $\frac{\lambda_{1,h}+\lambda^c_{1,h}}{2}$ in Figures 1 and 2. \\
\indent  From Table 1, on the one hand, we see that $\lambda_{1,h}$ converge to $\lambda_1$ from above and the corrected eigenvalue $\lambda^c_{1,h}$ converge to $\lambda_1$ from below, which indicate that the correction (\ref{s3.5}) provides lower eigenvalue bounds even though eigenfunctions are smooth. This coincides in the result of Theorem 3.1. On the other hand, on each domain, the CPU time to compute $\lambda^c_{1,h}$ is almost the same as that of $\lambda_{1,h}$, which tell us that the correction takes very little time. Furthermore, from Figures 1 and 2 we see that, on each domain, the error curves of $\lambda^c_{1,h}$ and $\lambda_{1,h}$ are parallel to the line with slope 2, which indicate $\lambda^c_{1,h}$ and $\lambda_{1,h}$ have the same and optimal convergence order $\mathcal{O}(h^2)$ and coincide in the result of Theorem 3.2. \\
\indent Although $\lambda^c_{1,h}$ are guaranteed to be the lower bounds of $\lambda_1$ whether the eigenfunctions are singular or smooth, compared with $\lambda_{1,h}$, the accuracy of $\lambda^c_{1,h}$ is slightly reduced. Thus, in order to  make up the loss of accuracy caused by correction, we use the average $\frac{\lambda_{1,h}+\lambda^c_{1,h}}{2}$ as a new approximation. From Figures 1 and 2 we know that, compared with $\lambda^c_{1,h}$,  $\frac{\lambda_{1,h}+\lambda^c_{1,h}}{2}$ have higher accuracy. Especially, for the square and the hexagon, the errors of $\frac{\lambda_{1,h}+\lambda^c_{1,h}}{2}$ are less than or equal to that of $\lambda_{1,h}$.
\begin{table}[h!]
  \centering
  \caption{The uncorrected eigenvalues and the corrected eigenvalues on $\Omega\subset \mathbb{R}^2$.}
    \begin{tabular}{ccccccc}
      \hline\noalign{\smallskip}
  domain &\multicolumn{2}{c}{\bf S}&\multicolumn{2}{c}{\bf L}&\multicolumn{2}{c}{\bf H}\\
$h$&$\lambda_{1,h}$ &$\lambda^c_{1,h}$&$\lambda_{1,h}$&$\lambda^c_{1,h}$&$\lambda_{1,h}$&$\lambda^c_{1,h}$   \\
\hline
$\frac{h_0}{32}$ &0.24008533  & 0.24006902  & 0.34143156  & 0.34134357  & 0.39334226  & 0.39329159  \\
$\frac{h_0}{64}$ &0.24008065  & 0.24007657  & 0.34141986  & 0.34139787  & 0.39332055  & 0.39330788  \\
$\frac{h_0}{128}$ &0.24007948  & 0.24007846  & 0.34141699  & 0.34141149  & 0.39331513  & 0.39331196  \\
$\frac{h_0}{256}$ &0.24007918  & 0.24007893  & 0.34141628  & 0.34141490  & 0.39331377  & 0.39331298  \\
$\frac{h_0}{512}$ &0.24007911  & 0.24007905  & 0.34141610  & 0.34141576  & 0.39331344  & 0.39331324  \\
$t(s)$ &31.10  &     31.20  &     22.74  &     22.81  &     25.34  &     25.41  \\
Trend&$\searrow$&$\nearrow$&$\searrow$&$\nearrow$&$\searrow$&$\nearrow$\\
   \hline
   \end{tabular}%
\end{table}%
\begin{figure}[h!]
\includegraphics[width=6cm,height=5cm]{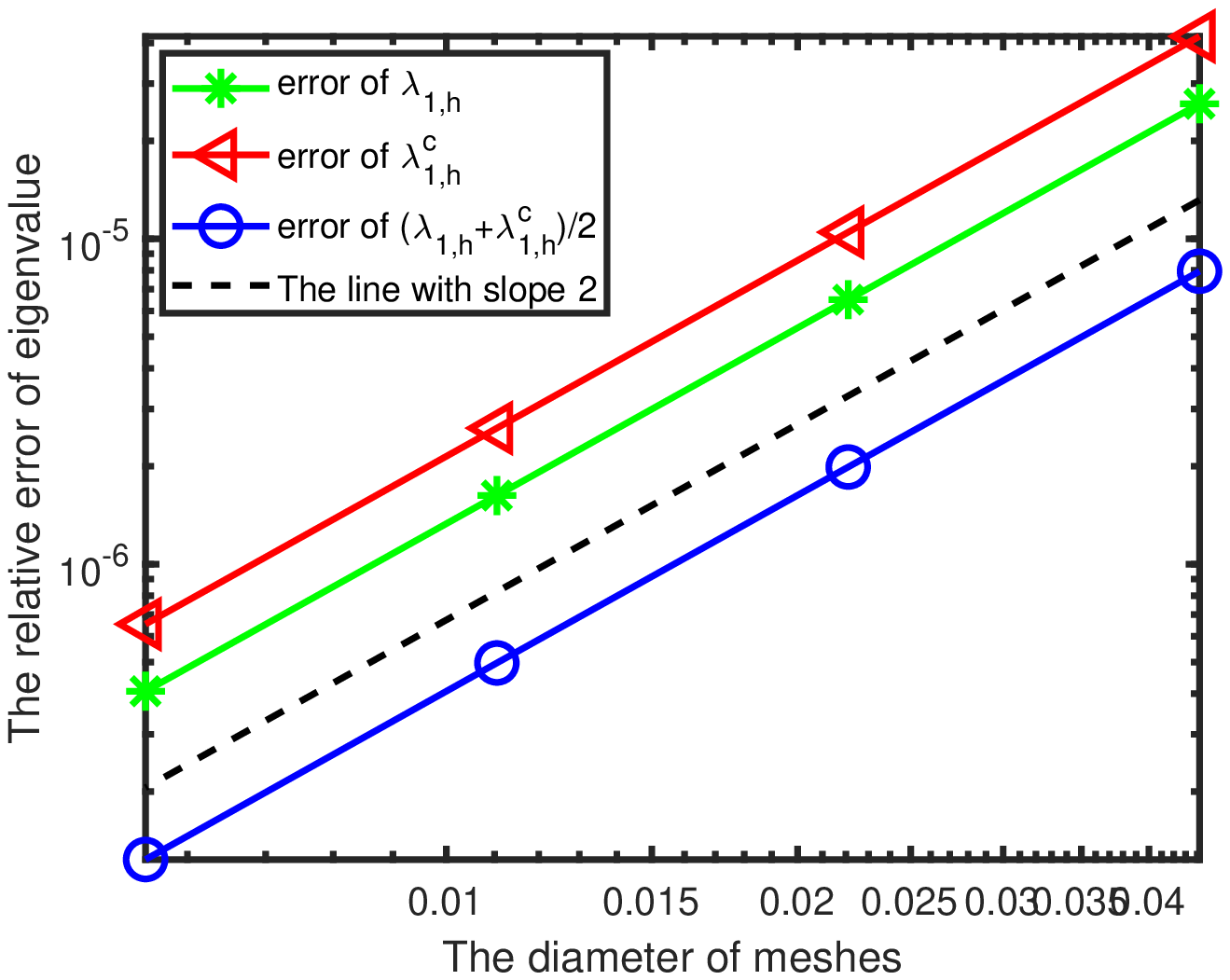}
\includegraphics[width=6cm,height=5cm]{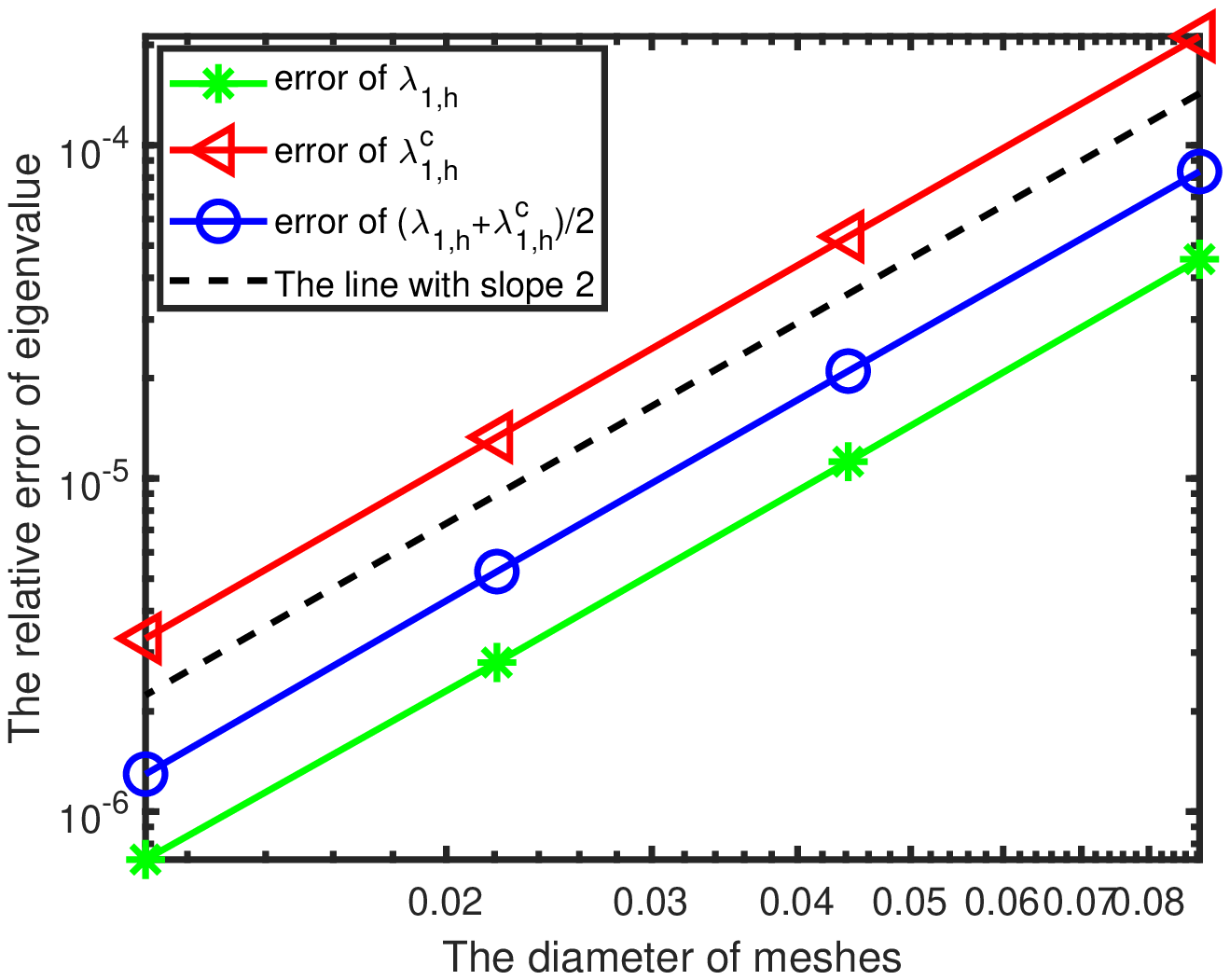}
\caption{\emph{The error curves of the first eigenvalues on the unit square (left) and the L-shaped domain (right) }}
 \end{figure}
 \begin{figure}[h!]
\includegraphics[width=6cm,height=5cm]{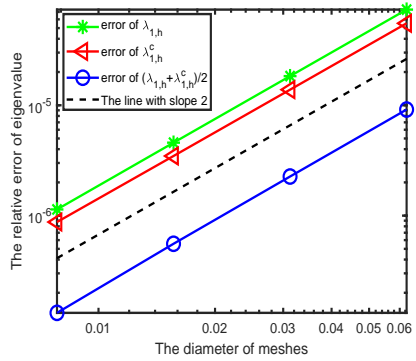}
\caption{\emph{The error curves of the first eigenvalues on the regular hexagon}}
 \end{figure}
\subsection{Numerical results on $\Omega\subset \mathbb{R}^3$}
\indent When $\Omega\subset \mathbb{R}^3$, we compute in the cube $(0,1)^3$ and the Fichera corner domain $[-1,1]^3\setminus(-1,)]^3$. For convenience, we simplify the domains as {\bf C} and {\bf F}, respectively. The quasi-uniform mesh samples of the cube and the Fichera corner domain are depicted in Figure 3. In the two domains, we compute the first three eigenvalues using CR finite element and list the results in Table 2. In the cube, $\lambda_2$ and $\lambda_5$ are eigenvalues with a multiplicity of 3. Corrected eigenvalues $\lambda^c_{1,h}$ are listed in Table 3. The error curves are depicted in Figure 4.\\
\indent  From Table 2, we see that, on each domain, $\lambda_{1,h}$ converge to $\lambda_1$ from above. From Figure 4, we know that the eigenfunctions corresponding to $\lambda_1$ are smooth. This shows that the CR finite element eigenvalue approximations are not necessarily lower bounds in the case of smooth eigenfunctions. From Table 3, we see that, corrected eigenvalues $\lambda^c_{1,h}$ converge to $\lambda_1$, which indicate that the correction (\ref{s3.5}) provides lower bounds for eigenvalues even though the eigenfunctions are smooth. From Figure 4, we see that the error curves of $\lambda^c_{1,h}$ and $\lambda_{1,h}$ are parallel to the line with slope 2, which indicate $\lambda^c_{1,h}$ and $\lambda_{1,h}$ have the same and optimal convergence order $\mathcal{O}(h^2)$.  The numerical results on three dimensional domains coincide in the result of Theorem 3.1 and Theorem 3.2.
\begin{figure}[h!]
\includegraphics[width=6cm,height=5cm]{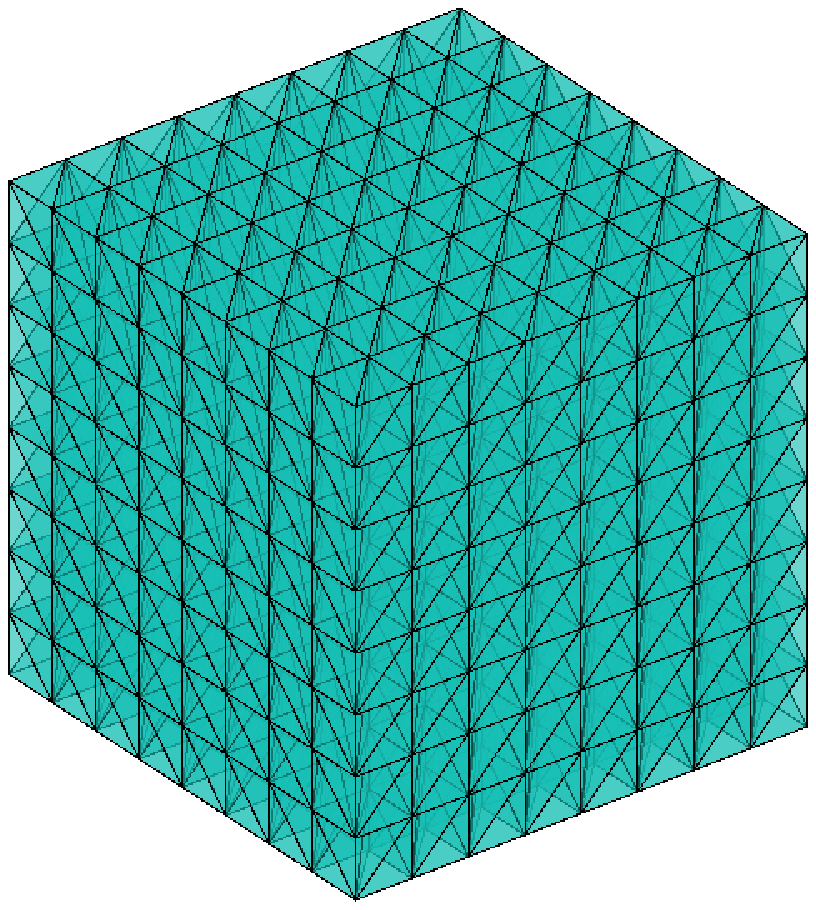}
\includegraphics[width=6cm,height=5cm]{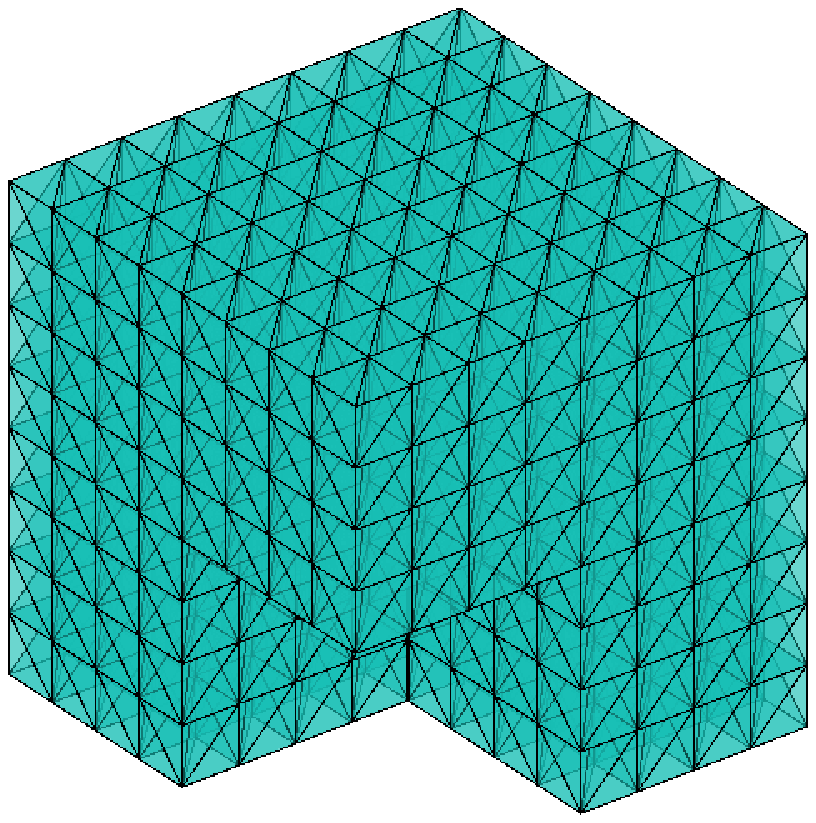}
\caption{\emph{The quasi-uniform mesh samples of the cube (left) and the Fichera corner domain (right) }}
 \end{figure}
 \begin{figure}[h!]
\includegraphics[width=6cm,height=5cm]{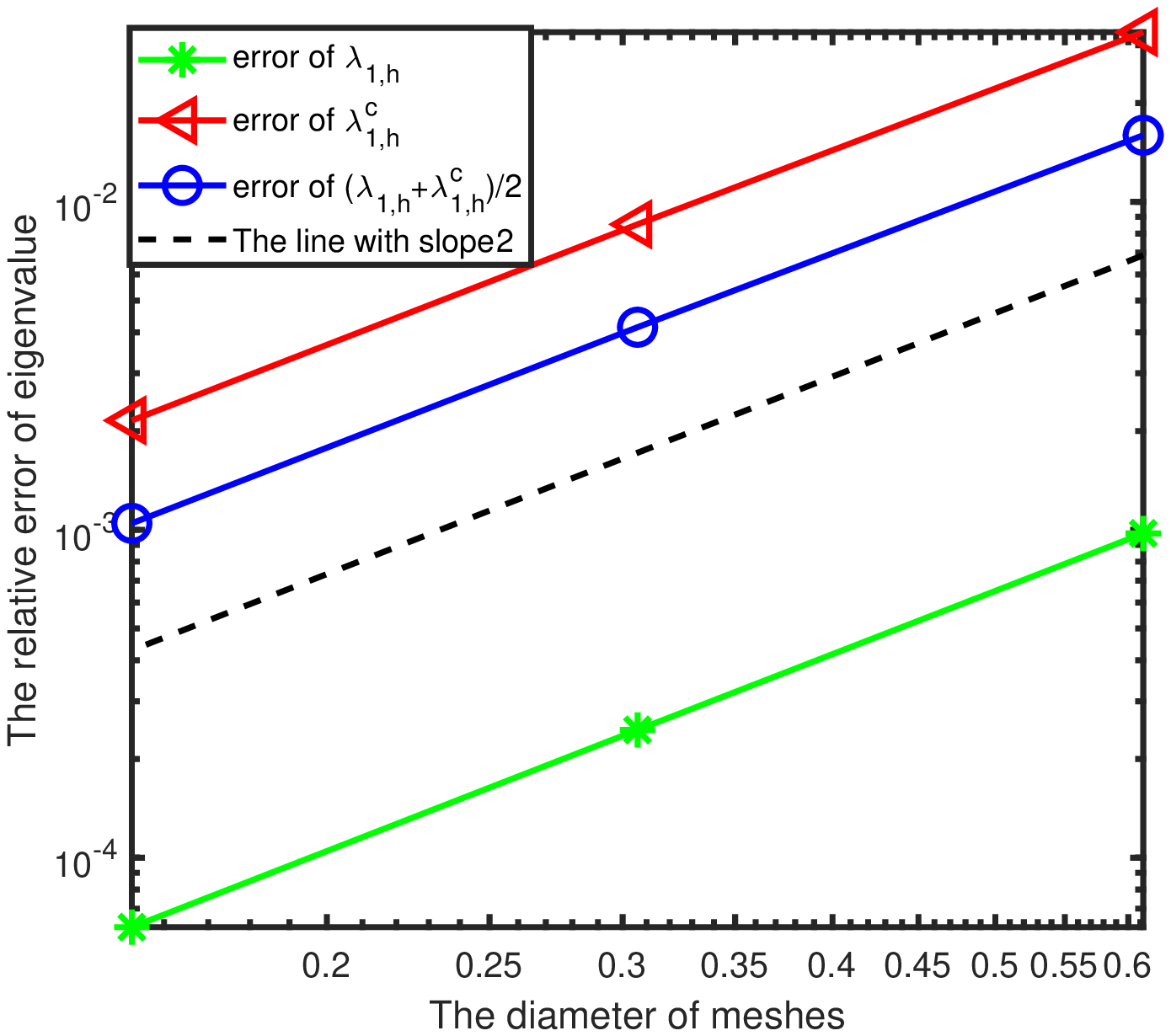}
\includegraphics[width=6cm,height=5cm]{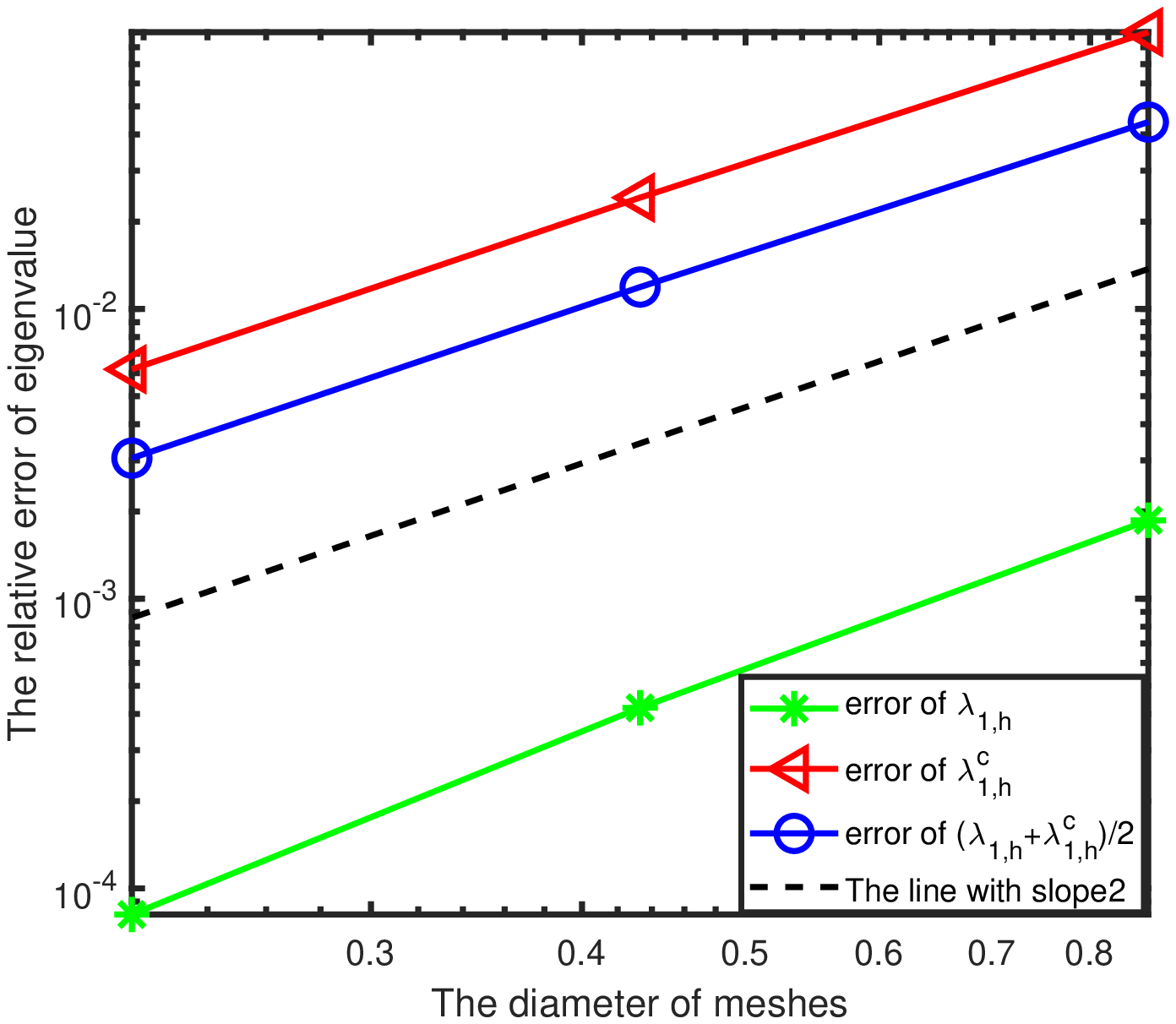}
\caption{\emph{The error curves of the first eigenvalues in the cube (left) and the Fichera corner domain (right) }}
 \end{figure}
  \begin{table}[h!]
  \centering
  \caption{The CR finite element eigenvalues approximations on $\Omega\subset \mathbb{R}^3$.}
    \begin{tabular}{ccccccccc}
      \hline\noalign{\smallskip}
  domain &\multicolumn{3}{c}{\bf C}&\multicolumn{4}{c}{\bf F}\\
$h$&$\lambda_{1,h}$ &$\lambda_{2,h}$&$\lambda_{5,h}$&$h$&$\lambda_{1,h}$ &$\lambda_{2,h}$&$\lambda_{3,h}$  \\
\hline
0.6124  &  0.1623444 &   1.11356  &   1.56489  &    0.8660  &   0.268747 &   0.54947  &   0.72763  \\
0.3062  &  0.1622261 &   1.14537  &   1.65619  &    0.4330  &   0.268359 &   0.56641  &   0.73377  \\
0.1531  &  0.1621963 &   1.15272  &   1.68222  &    0.2165  &   0.268268 &   0.57235  &   0.73615  \\
0.0765  &  0.1621888 &   1.15448  &   1.68924  &    0.1083  &   0.268247 &   0.57441  &   0.73687  \\
Trend&$\searrow$&$\nearrow$&$\nearrow$&--&$\searrow$&$\nearrow$&$\nearrow$\\
   \hline
   \end{tabular}%
\end{table}%
 \begin{table}[h!]
  \centering
  \caption{The uncorrected eigenvalues and the corrected eigenvalues on $\Omega\subset \mathbb{R}^3$.}
    \begin{tabular}{cccccc}
      \hline\noalign{\smallskip}
  domain &\multicolumn{2}{c}{\bf C}&\multicolumn{3}{c}{\bf F}\\
$h$&$\lambda_{1,h}$ &$\lambda^c_{1,h}$&h&$\lambda_{1,h}$&$\lambda^c_{1,h}$   \\
\hline
0.6124  &   0.162344 &   0.156854 &    0.8660  &   0.268747 &   0.244062 \\
0.3062  &   0.162226 &   0.160802 &    0.4330  &   0.268359 &   0.261752 \\
0.1531  &   0.162196 &   0.161837 &    0.2165  &   0.268268 &   0.266587 \\
0.0765  &   0.162189 &   0.162099 &    0.1083  &   0.268247 &   0.267824 \\
$t(s)$&    150.07  &    150.23  &      -- &    223.12  &    223.26  \\
Trend&$\searrow$&$\nearrow$&--&$\searrow$&$\nearrow$\\
   \hline
   \end{tabular}%
\end{table}%

{\small
{\em Authors' addresses}:\\
{\em Yu Zhang}, School of Mathematical Sciences, Guizhou Normal University, Guiyang 550001, China;\\School of Mathematics $\&$ Statistics, Guizhou University of Finance and Economics, Guiyang 550001, China.
 e-mail: \texttt{zhang\_hello\_hi@126.com}.\\
{\em Hai Bi}, School of Mathematical Sciences, Guizhou Normal University, Guiyang 550001, China.
 e-mail:  \texttt{bihaimath@gznu.edu.cn}.\\
{\em Yidu Yang (corresponding author)}, School of Mathematical Sciences, Guizhou Normal University, Guiyang 550001, China.
 e-mail:  \texttt{ydyang@gznu.edu.cn}.
}

\end{document}